\documentstyle[11pt]{article}

\input{psfig}

\begin{document}
\thispagestyle{empty}

\def\qed{\nobreak\hskip 2pt\vrule height 5pt width 5pt depth 0pt}
\def\proof{\noindent{\sl Proof: }}
\def\Remark{\noindent {\bf Remark:  }}
 \def\D{{\bf D}}
 \def\R{{\bf R}}
  \def\T{{\bf T}}
  \def\Z{{\bf Z}}
 \def\Q{{\bf Q}}
 \def\N{{\bf N}}
 \def\C{{\bf C}}
 \def\TTT{{\cal T}}
 \def\qed{{\bf Q.E.D.}}
 \def\O{{\bf O}}
 \def\DD{{{I\negthinspace\!D}}}
 \def\RR{{{I\negthinspace\!R}}}
 \def\NN{{{I\negthinspace\!N}}}
 \def\TT{{{T\negthinspace\negthinspace\!T}}}
 \def\ZZ{{{Z\negthinspace\negthinspace\!Z}}}
 \def\QQ{{{0\negthinspace\negthinspace\negthinspace\!Q}}}
 \def\CC{{{C\negthinspace\negthinspace\negthinspace\!C}}}
 \def\SS{{{\tilde S}}}
 \def\Cstar{{\CC^*}}
\def\Chat{{\hat{{C\negthinspace\negthinspace\negthinspace\!C}}}}

 \begin{titlepage}
 \title{ The Set of Maps $F_{a,b}:x \mapsto x+a+{b\over 2\pi}
   \sin(2\pi x)$ with any
 Given  Rotation Interval is Contractible}

 \author{
Adam Epstein\\Mathematics Department\\California Institute of Technology\\
Pasadena, CA 91125\and
 Linda Keen\thanks{ Supported in part by NSF GRANT
DMS-9205433, Inst. Math. Sciences, SUNY-StonyBrook and I.B.M.}
\\Mathematics Department\\CUNY Lehman College\\Bronx, NY 10468 \and
 Charles
Tresser\\I.B.M.\\ PO Box 218\\ Yorktown Heights N.Y. 10598}

 \maketitle
 \begin{abstract}
Consider the two-parameter family of real
analytic maps $F_{a,b}:x \mapsto x+ a+{b\over 2\pi} \sin(2\pi x)$ which are
lifts of degree one endomorphisms of the circle. The purpose of this paper
is to provide a proof that for any closed interval $I$, the set of maps
$F_{a,b}$ whose rotation interval is $I$, form a contractible set.
\end{abstract}
\end{titlepage}
 \bibliographystyle{plain}

\section{Introduction}
\label{sec:intro}

  Orientation preserving homeomorphisms and diffeomorphisms of the
circle have attracted the attention of mathematicians and physicists
for a long time because they arise as Poincar{\'e} maps induced by
non-singular flows on the two-dimensional torus \cite{Ar,De,Herman,Po}. More
recently, families of circle endomorphisms which are deformations of
rotations have appeared as approximate models for some scenarios of
transition to ``chaos", or more technically, transitions from zero to
positive topological entropy. One can observe these transitions by
varying parameters of flows in $3$-space so that tori supporting
non-singular flows get wrinkled and then get destroyed. More generally,
these scenarios are typical for a huge variety of systems of coupled
oscillators so that one sees them everywhere. As a matter of fact, in
many cases when one is lead to study an endomorphism of the interval as
a model for a  natural science experiment, some circle endomorphism is
a more adequate model.

  While the simplest endomorphisms of the
interval depend on a single parameter, say the non-linearity, the
simplest reasonably complete family of circle endomorphisms containing
the rotations has to depend on two parameters: the non-linearity and
some form of mean rotation speed.  In the coupled oscillators picture,
these parameters correspond respectively to the strength of the
forcing and the frequency ratio of the coupled oscillators. The
paradigm for interval endomorphisms is the quadratic family. For the
circle, it is the so-called {\em Arnold} or {\em standard} two-parameter
family:

 $$f_{A,b}:\theta \mapsto (\theta +A+{b\over 2\pi}
\sin(2\pi\theta))_1\,,$$
 with $(A,\,b)\in [0,\,1[\times \RR^+$, and
$(\theta)_n\buildrel \rm def\over = \theta\,{\mathop {\rm mod}}\,n.$

 Under an orientation preserving homeomorphism of the circle, the orbits of
all points wrap around the circle at the same average speed \cite{Po}. For
non-invertible maps this is no longer necessarily the case, but the
set of average speeds form a closed interval. With the concepts roughly
recalled so far, we can give an example of a result that is a corollary of
our main theorem: for any average speed $\omega$, the set $\lbrace (A,\,b)
\rbrace$ of pairs such that all orbits under the map $f_{A,b}$ wrap
around the circle at speed $\omega$, is connected. Our main result is
in fact a similar statement in the more general setting where average
speeds vary in an interval.

     Precise definitions and statements are contained in $\S 2$.  In
$\S 3$, we reduce our main theorem to a rigidity result: this
reduction is merely well known material, but some proofs are sketched
for completeness.  In $\S 3$ we have also included  some material not
strictly needed for the proof of Theorem A, but intended to help some
readers to build an intuitive picture of what the main result is all
about. The rigidity property, formalized in Theorem D, is  proved in
$\S 4$.  Our proof of Theorem D is
one more example of the efficiency of complex analytic methods in
dealing with questions arising naturally in a real analytic framework.

\section{ Definitions and statement of the results}
\label{sec:defns}

  Let $\TT=\RR /\ZZ$ be the circle and $\Pi:\RR\to \TT$ the canonical
projection.  The real continuous map $F$ is a lift of the continuous
circle map $f:\TT\to \TT$ if and only if
$$f\circ\Pi=\Pi\circ F\,.$$
The integer $d$ such that $$F(x+1)=F(x)+d\,\,,$$
for all real numbers $x$ is called the {\it degree} of $f$ (or of
$F$). The identity map, and more generally the rotations, have degree
one. Since the degree varies continuously for continuous deformations
of circle maps, and since we are interested in a parametrized continuous
family containing rotations, we shall only consider degree one maps in
the rest of the paper.  Hence,  {\it circle map} will always
mean ``degree-one continuous circle map", and a real map
will be called a {\it lift} if and only if it is the lift of a degree
one circle map.

  Let $f$ be a circle map, and let $F$ be a lift of $f$
(each time both symbols $f$ and $F$ appear conjoined in the paper,
they are related in the same way). We define
$$\underline{\rho}_F(x)=\liminf_{n\to\infty}{F^n(x) \over n}\,\,,$$
and
$$\overline{\rho}_F(x)=\limsup_{n\to\infty}{F^n(x) \over n}\,\,.$$
The {\it rotation interval} of $F$ \cite{NPT} is then
$$I(F)=[\,\alpha,\,\beta]\,\,,$$
where
 $$\alpha=\inf_{x\in
\RR}\underline{\rho}_F(x)\,\,,\,\,\beta=\sup_{x\in
\RR}\overline{\rho}_F(x)\,\,.$$ When $I(F)$ is a singleton $\lbrace
\omega\rbrace$, we sometimes use the classical language and say that
$\rho (F)\buildrel \rm def\over =\lbrace \omega \rbrace$ is the {\em
rotation number} of $F$.

We will focus on the standard family $f_{A,b}$ with parameter space
$[0,1[\times\RR ^+,$ and the corresponding degree one lifts $F_{a,b}$ with
parameter space $\RR\times\RR ^+,$ where the correspondence is given by $A=a\,
\,{\mathop{\rm mod}}\,\,1$. To state  our main result we need the following
 \proclaim Definition. An arc or curve $a =
\phi(b)$ in parameter space $\RR \times \RR^+$ is called an {\em L-curve} if
$\phi$ is uniformly Lipschitz  with  bound $1\over 2 \pi$.

 \bigskip
 \proclaim Theorem A.  For each closed interval $I$, the set
$R_I$ of standard lifts with rotation interval $I$ corresponds to a
contractible region, also denoted $R_I$, in the parameter space
$\RR\times \RR^+$. More precisely,
\begin{itemize}
 \item For any
irrational number $\omega$, $R_{\lbrace \omega\rbrace}$ is an L-curve.
 \item If one bound of
$I$ is irrational while the second bound is a rational number, $R_I$
is an L-curve.
 \item If
the bounds of $I$ are distinct irrational numbers, $R_I$ is an L-curve.
 \item When both bounds of $I$ are rational, $R_I$ is a lens shaped domain
bounded  by two
L-curves that meet at their endpoints.
 \end{itemize}

 \bigskip

\noindent {\bf Convention}.
To simplify the language and the notation, we shall continue to
identify sets of standard lifts with the corresponding regions in
parameter space, as we did in the statement of Theorem A; when the distinction
is relevant the context should tell which space we mean.

\bigskip
\proclaim Conjecture B.  If the bounds of $I$ are distinct
irrational numbers, $R_I$ is a point.

 \bigskip
The content of Theorem
A is illustrated in Figures 1 to 3. The labeling appearing in these Figures
is defined in $\S\S$ 2 and 3.

\begin{figure}[htp]
  \centerline{\psfig{figure=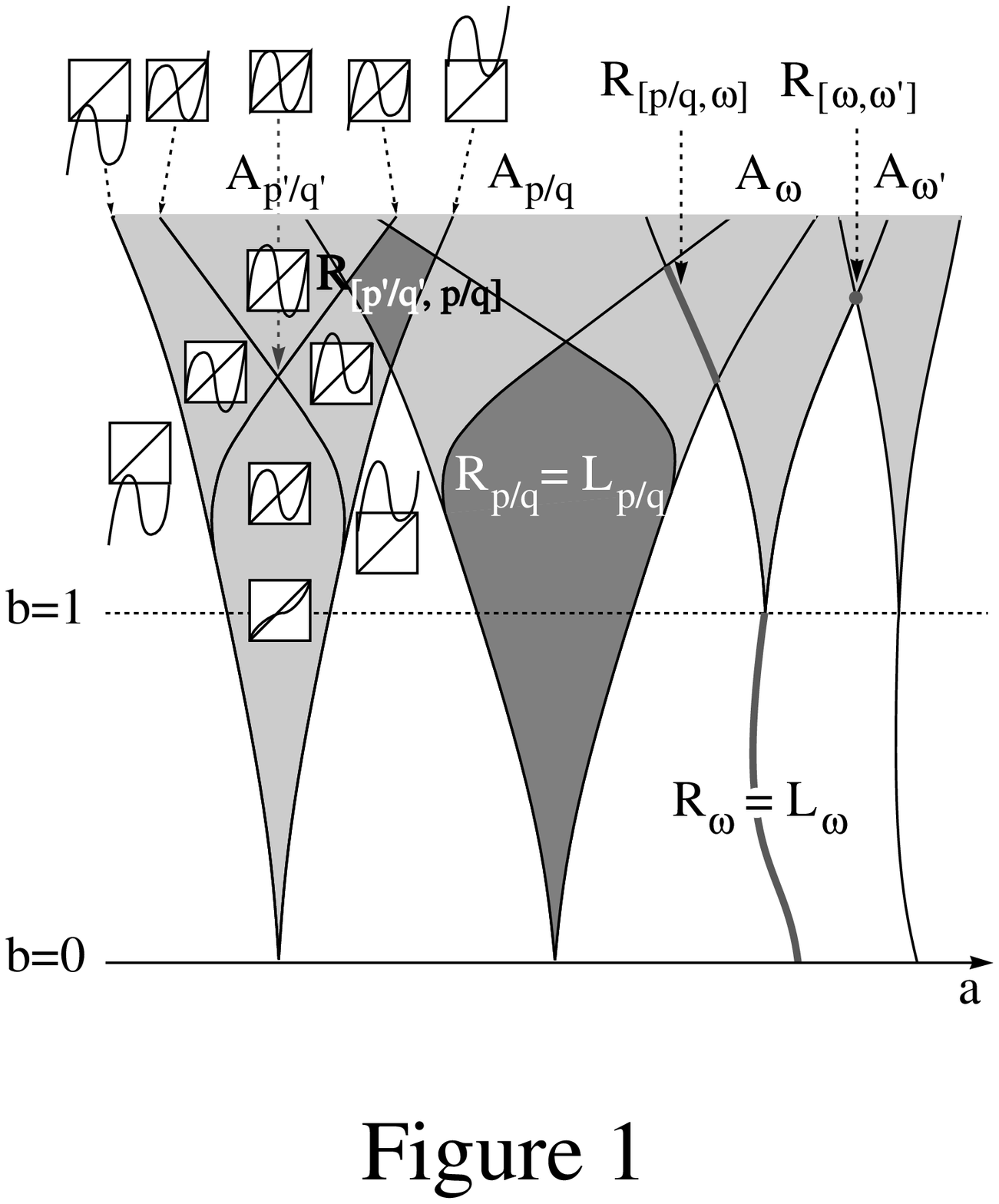,width=.4\vsize}}
  \smallskip
   \begin{quote}
     A schematic picture of part of the parameter space for the lifts
     of the standard family. The small inserts represent the graph of $F^{q'}$
     for $F$ in  various regions,  lines and points.
   \end{quote}
  \bigskip
  \centerline{\psfig{figure=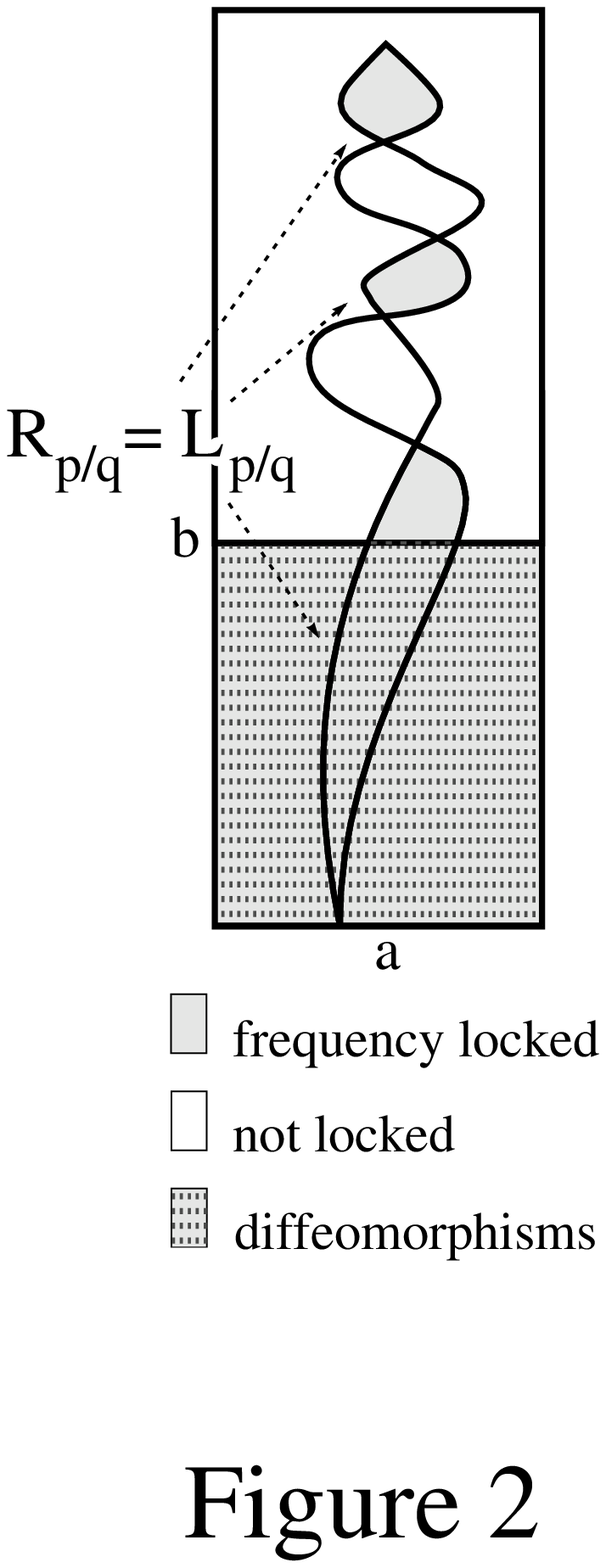,height=.3\vsize} \hfil
              \psfig{figure=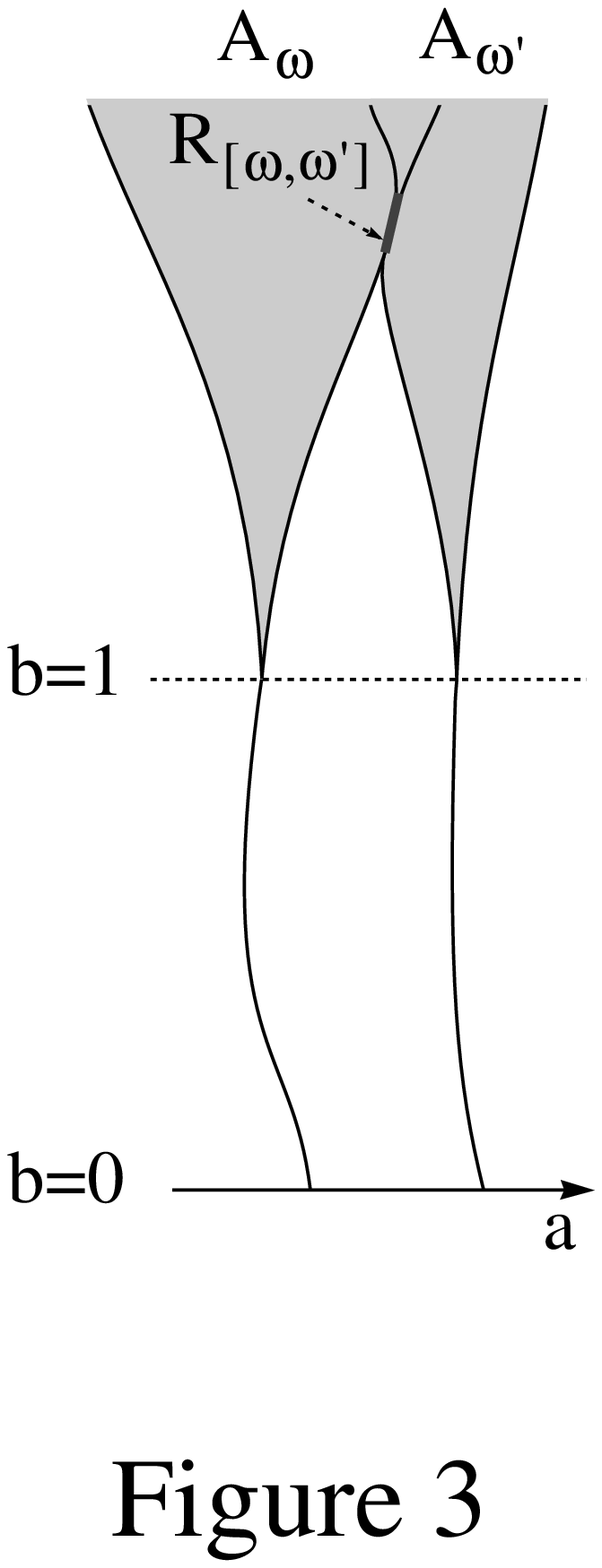,height=.3\vsize}}
  \smallskip
  \begin{quote}
  	Figure~2 shows a situation proved to not exist in the standard
         family by Theorem~A, and 
         Figure~3 illustrates the situation conjectured not to occur in
         the standard family (in Conjecture~B).
  \end{quote}
\end{figure}

 As an important particular case, relevant for example in the description of
the boundary of chaos \cite{MaT}, Theorem A contains the following result which
was conjectured in \cite{Bo} (p. 378 and Fig. 13) and \cite{MaT} (p. 213):
 \bigskip
\proclaim Corollary C.  For each real $\omega$, the set $R_{\lbrace
\omega\rbrace}$ is connected.

 \bigskip
\Remark
  Corollary C
(which was known for a long time to hold in the case when $\omega$ is
irrational) was recently proved for some families of piecewise affine
lifts; for theses families, the proof only uses elementary real
variable methods \cite{UTGC}. Later, and in parallel to the work presented
here, more sophisticated real variable methods \cite{GMT} were used to get the
counterpart to Theorem A as well as a proof of Conjecture B for these piecewise
lifts.

\section{ Proof of Theorem A, Part I:  Real analytic part}
\label{sec:real analysis}

\subsection{ Non-decreasing lifts}

  We shall denote by ${\cal F}^k(\RR)$ the space of $C^k$ non-decreasing
lifts, equiped with the $C^k$ norm. The next theorem recalls some classical
results \cite{Ar,Herman}, usually formulated for lifts of degree one
homeomorphisms, but generalizable {\it verbatim} to ${\cal F}^0(\RR)$.

\bigskip

\proclaim Theorem 3.1.1.
\begin{enumerate}
\item
 The rotation number, as a function $\rho:{\cal
F}^0(\RR)\to \RR$
  is continuous.\\
   \item
  For $F$ and $G$ in
${\cal F}^0(\RR)$,
  $$F\geq G\qquad\Rightarrow\qquad \rho
(F)\geq
\rho(G)\,\,,$$
  $$F>G\,\,\&\,\,\rho(F)\,\mbox{or\,}\,\rho(G)\,\mbox{irrational}
\quad\Rightarrow\quad
\rho (F)>
  \rho(G)\,\,.$$
  \item
 If $\rho(F)$ is irrational
and
  $f$ has a dense orbit, then, $$F\geq G\quad and\quad
  F\not=
G\qquad\Rightarrow\qquad \rho (F)>\rho(G)\,\,.$$
\end{enumerate}

  \medskip
  The next
result which can also be found in the above references, is more
specific
  to our problem:
\proclaim Theorem 3.1.2.
For $b\not = 0$, no iterate of a standard lift is
affine.

  \medskip
  Set
  $${\bf A}'_\omega=\lbrace F\in {\cal
F}^0(\RR)\,|\,\omega\in I(F) \rbrace\,,$$

  \medskip

  Theorems 3.1.1 and 3.1.2 yield a
partition of the subset
  $\RR\times [0,1]$ of the parameter space of
the standard family in the ${\bf
  A}'_\omega$'s, with the following
properties \cite{Ar}:

 P1. For $\omega$ irrational, ${\bf
A}'_\omega$ is an arc crossing each line
  $b=constant\leq 1$ at a
single point,

P2. For a rational number ${p\over q}$,
${\bf A}'_{p\over q}$, is often called an
  {\it Arnold tongue}; it
crosses each line $0<b=\mbox{constant} \leq 1$ on an
  interval of positive length.

  \medskip

\Remark Property P2
 describes an aspect of the phenomenon of
  ``frequency locking'',
first described by Huyggens, in the context of clocks hanging
  from
the same wall, and descibed in modern terms, in the simplest cases, as
the
  strucural stability of generic degree one circle diffeomorphisms
with rational
  rotation numbers.

\subsection {Some special sets}

Following \cite{Bo} and \cite{MaT}
(both of whom extended the above mentioned work in \cite{Ar} from
homeomorphisms to endomorphisms), for each real number $\omega$, we define

  $${\bf A}_\omega=\lbrace F \in
  \lbrace F_{a,b}\rbrace_{(a,\,b)\in
\RR\times
  \RR^+}\,\,|\,\omega\in I(F) \rbrace\,,$$
  and
  $${\bf L}_\omega=
\lbrace F \in \lbrace F_{a,b}\rbrace_{(a,\,b)\in
\RR\times
  \RR^+}\,\,|\,\lbrace \omega \rbrace
=I(F)
  \rbrace\,.$$

  According to our notational convention, ${\bf A}_\omega$ and ${\bf
L}_\omega$ can also be understood as subspaces of the parameter space. The
following theorem enables us to use the results of \S~3.1 to analyze these
subspaces.

\medskip
\proclaim  Theorem 3.2.1 {\rm (\cite{CGT,Mi})}.
  \begin{enumerate}
\item
 For any lift $F\in
C^0_1(\RR)$,
  $$I(F)=[\,\rho(F^-),\,\rho(F^+)]\,,$$
  where
  $F^+$ is
the monotonic upper-bound of $F$, and $F^-$ is the
monotonic
  lower-bound (see Figure 4). In formulas we
have:
  $$F^+(x)=sup_{y\leq x}(F(y))\,,$$
  $$F^-(x)=inf_{y\geq
x}(F(y))\,.$$
\item
 For each $\omega\in I(F)$, there
is a non-decreasing
  lift $F_\omega$ with $\rho(F_\omega)=\omega$, and
such that $F_\omega$ coincides
  with $F$ where it is not locally
constant.
\end{enumerate}

\begin{figure}[ht]
	\centerline{\psfig{figure=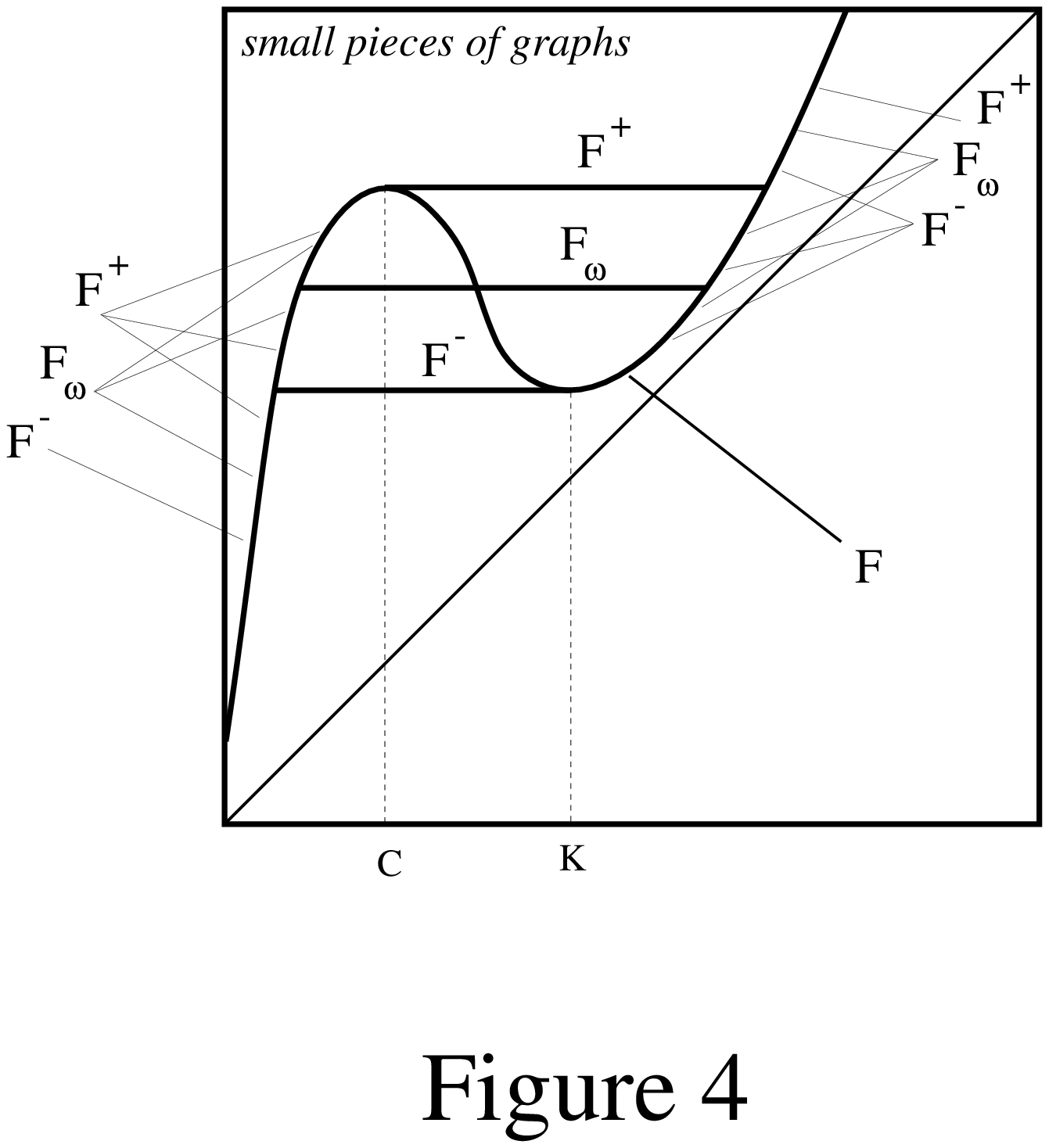,width=.5\hsize}}
\end{figure}

  \medskip
  Theorem 3.2.1. is quite easy to prove for the
 maps in the standard family.

\subsection{ Simple properties of the ${\bf
A}_\omega$'s and ${\bf
  L}_\omega$'s.}

  Theorem 3.2.1. allows us to
use Theorem 3.1.1. in the study of non-invertible
  maps.  In
particular, one gets easily that the ${\bf A}_\omega$'s are
connected
  and, for $\omega$ rational, they intersect each
line
  $b=\mbox{constant} >1$ on a segment of non-zero length. We discuss the
case when $\omega$
 is irrational in \S~3.4. A fundamental role
will be played by the boundaries of
  the ${\bf A}_\omega$'s and some of
their accumulation sets defined as follows:
 $${\bf B}_{\omega}^l=\lim _{{\theta}\to{\omega}^+}{\bf A}_{\theta}^l\,,$$
 $${\bf B}_{\omega}^r=\lim _{{\theta}\to{\omega}^-}{\bf A}_{\theta}^r\,.$$
 These boundaries, pieces of which form all the boundaries of the sets
$R_I$ of Theorem A, are
  described in the
following  result due to Boyland \cite{Bo}. We include a proof for the sake of
completeness.

\medskip
\proclaim Theorem 3.3.1(\cite{Bo}).

\begin{enumerate}
\item
 For any
real number $\omega$, the left and right bounds ${\bf
A}_{\omega}^l$
  and ${\bf A}_{\omega}^r$ of ${\bf A}_{\omega}$ in
$\RR\times \RR ^+$, are L-curves.
\item
When $\omega$ is a
rational number these curves intersect at a
single point which has $b=0$ as second coordinate
\item

  $${\bf B}_{p\over q}^l=\lim _{{p'\over q'}\to{p\over q}^+}{\bf
A}_{p'\over q'}^l\,,\, \mbox{\rm  and }{\bf B}_{p\over q}^r=\lim _{{p'\over
q'}\to{p\over q}^-}{\bf A}_{p'\over q'}^r\,.$$

Moreover, the sets ${\bf B}_{p\over q}^l$ and ${\bf B}_{p\over q}^r$ are
L-curves.
\item
 For $\omega$ irrational
  $${\bf
A}_\omega ^l={\bf B}_\omega ^l=\lim _{{p\over q}\to\omega ^-}
{\bf
  A}_{p\over q}^l\,,$$
  $${\bf A}_\omega ^r={\bf B}_\omega
^r=\lim _{{p\over q}\to\omega ^+} {\bf
  A}_{p\over q}^r\,.$$
\end{enumerate}

  \bigskip

\proof
With no loss of generality, we prove statement~1 for ${\bf A}_{p\over q}^l$.
To do this we consider the vertical cone in $\RR\times \RR ^+$ with vertex at
$(a,b)$ and boundaries made by lines with slopes $2\pi$ and $-2\pi$, and show
that it contains ${\bf A}_{p\over q}^l$. A point to the right of the cone
has coordinates $a'=a+\delta+\epsilon$, $b'=b\pm 2\pi \delta$, while a point to
the left of the cone has coordinates $a"=a-\delta-\epsilon$, $b"=b\mp 2\pi
\delta$, for some non-negative $\delta$ and $\epsilon$. Consequently, using the
continuity of the rotation number of $F_{a,b}^+$ as a function of the
parameters, and its monotonicity as a function of $a$, the inclusion  of ${\bf
A}_{p\over q}^l$ in the cone  follows from  $$\forall
\delta \geq 0,\,\, \forall \epsilon\geq 0:\quad \delta+\epsilon\pm \delta \sin
2\pi x\geq 0\,.$$
  Statement~2
 follows from from Theorem 3.1.2. Statement~3
 follows from statement~1
by continuity of the rotation number applied to the
  monotonic bounds.
Statement~4 is a consequence of the same
continuity
  property.
  \qed

\subsection{Theorem A in the simplest case}

  We recall here the proof of Theorem A in the case when
$I$ is the singleton
  $\lbrace \omega\rbrace$ for some irrational
number $\omega$. We begin with a weak
  form of a theorem by Denjoy
\cite{De}

\proclaim Theorem 3.4.1 (\cite{De}).
  For $F$ in ${\cal F}^2(\RR)$,
  if $\rho(F)= \omega $ for
some irrational number $\omega$, then the
  circle map $f$ with lift
$F$ has a dense orbit.

  \medskip
  The next result is a particular
case of a theorem obtained by Block and
  Franke (see also \cite{CGT}) as a
consequence of the Denjoy theory:

\proclaim Theorem 3.4.2 (\cite{BF}). If $b>1$ and $\rho = \rho
(F_{a,b}^-)=\rho (F_{a,b}^+),$ then $\rho \in \QQ$.

 \proof  We
  first remark that there exist distinct
$C^2$ smooth lifts $F_0$ and $F_1$ such
  that,
$\forall
  x\,\in\,\RR,\,\,\,F^-(x)\,\leq\,F_0(x)\,\leq\,F_1(x)\,\leq\,F^+(x)$.
If the
  claim were false, by Theorem 3.1.1-2,
$\rho=\rho(F_0)\,=\,\rho(F_1)\,=\,\rho(F^+)$ for some $\rho \notin
  \,\QQ$. But if $\rho \not\in \QQ$,
Theorem 3.4.1 implies $F_0$ (and $F_1$) has a dense orbit. Hence the
  claim
follows from Theorem 3.1.1-3. \qed

  \medskip
  To finish the analysis of the case when $\rho(F)= \omega\not\in \QQ$, we
just have
 to check that, as a consequence of Theorem 3.4.2, all
${\bf L}_{\omega}$ are contained in the region $b\leq1$ described in \S~3.1 and
\S~3.3. In summary, we have: \medskip

\proclaim  Lemma 3.4.3. For $\omega\not\in \,\QQ$, $R_{\lbrace \omega\rbrace}$
  is an L-curve
 contained in the
region $b\leq 1$.

\subsection{Intersections of the boundaries of the $A_{\omega}$'s:
existence}

  For any $A$, the narrowest diagonal strip with sides parallel to and
centered on the main diagonal that  contains the graph of $F_{a,b}$, can
be made arbitrarily wide by choosing $b$ large enough. Hence
 \medskip

\proclaim  Lemma 3.5.1. For any $\omega\in\RR$, and any $a\in \RR$, $\omega$ is
contained in the interior of $I(F_{a,b})$ as soon as $b$ is large enough.

  \medskip
\proclaim Corollary 3.5.2.
If $\omega<\theta$, ${\bf A}_\omega ^r$ and ${\bf B}_\omega ^r$ intersect ${\bf
A}_\theta ^l$ and ${\bf B}_\theta^l$. Also, ${\bf B}_{p\over q}^r$ intersects
${\bf B}_{p\over q}^l$ for any rational ${p\over q}$.

\subsection{Intersections of the boundaries of the $A_{\omega}$'s:
combinatorics.}

In order to prove Theorem A using
Theorem 3.3.1-4,  we may restrict our attention to the intersection points
of the boundaries of $A_{p\over q}$ and $A_{p'\over q'}$. Before we begin the
analysis of these intersection points, let us recall that the Schwarzian
derivative of a map $g$ is defined as
  $$Sg={g'\negthinspace '\negthinspace'\over g'}-{3\over
2}({g'\negthinspace '\over g'})^2\,.$$
A direct computation then gives

  \medskip
\proclaim  Lemma 3.6.1. For $b>1$ and $F'_{a,b}(x)\neq 0$, $SF_{a,b}(x)<0$.

  \medskip
This will be used to prove Lemma 3.6.2.

Using a lift $O$ of any periodic orbit $o$ of $f_{A,b}$, we can
analyze the local behavior of $f_{A,b}$ near $o$ in terms of the derivatives of
$F_{a,b}$ at $q$ successive points of $O$, $x_0,x_1=
F_{a,b}(x_0),\ldots,x_{q-1}=F_{a,b}(x_{q-2})$. Define the {\em
multiplier} of $o$ as $m_o=F'_{a,b}(x_0)\cdot F'_{a,b}(x_1)\cdots
F'_{a,b}(x_{q-1})$. We call $o$

\noindent
- {\em attracting} if $|m_o|<1$,

\noindent
- {\em neutral} if  $|m_o|=1$, and in particular {\em parabolic} if $m_o$ is a
root of unity,

\noindent
- {\em hyperbolic} if $|m_o|>1$.

Clearly $m_o$ only depends on $O$.

The periodic orbits of the circle map $f_{A,b}$, which are also orbits of
some homeomorphism of the circle, and lifts of these orbits will play an
important role in our discussion. If a point of such a periodic orbit $o$ of
$f_{A,b}$ has a lift with rotation number ${p\over q}$ under $F_{a,b}$, $o$ has
period $q$ and lifts to $p$ distinct orbits of $F_{a,b}$.

\proclaim  Lemma 3.6.2. Suppose ${p\over q}\in I(F_{a,b})$. Then $F_{a,b}$ has
an orbit $O$ such that
\begin{enumerate}
\item
$O$ projects to a periodic orbit $o$ of $f_{A,b}$,
\item
there is a monotone $G\in {\cal F}^0(\RR)$ such that $O$ is an orbit of $G$,
\item
no point of $O$ is in an interval where $F_{a,b}$ is decreasing.
\item
If $m_o\geq 1$, $O$ is uniquely determined up to integer translation. If we
relax the multiplier condition, there are, up to integer
translation, at most two distinct orbits. When there are two
orbits distinct
 under integer translation, denote them by $O$ and $O'$. Then
\begin{enumerate}
\item
$O$ and $O'$ bound intervals that are lifts of $q$ pairwise disjoint arcs
on which $f_{A,b}$ is orientation preserving,
\item
the interiors of these intervals are in the immediate basin of the attracting
periodic orbit which lifts to $O'$,
\item
at least one critical point is in the immediate basin of the attracting
orbit
which lifts to $O'$.
\end{enumerate}
\end{enumerate}

  \medskip
\proof  For properties 1 to 3, the
existence follows from Theorem 3.2.1.

Uniqueness under the condition $m_o\geq 1$ follows from the well known fact
(proved by a direct computation) that the absolute value of the (usual)
derivative of a function on the real line, whose Schwarzian derivative is
negative off the critical set, has no local non-zero minimum. That there
are at most two orbits distinct under integer translation when the
multiplier condition is relaxed follows in a similar fashion from the
negative Schwarzian derivative property.

Property 4 (a) comes from the fact that we only use the restriction of
$F_{a,b}$ to the intervals where it is increasing; property 4 (b) is
immediate (draw a graph); and property 4 (c) is a classical result in
holomorphic dynamics (see the Remark in \S 4.1).  \qed

  \medskip
 Let
 $${\bf o}_{P\over Q}=\lbrace {\bf p}_0,\,{\bf p}_1,\dots,
{\bf p}_{Q-1}\rbrace\,\,\,{\rm and}\,\,\, {\bf o'}_{P\over Q}=
\lbrace{\bf p}'_0,\,{\bf p}'_1,\dots, {\bf p}'_{Q-1}\rbrace$$
 be the projections of
${\bf O}_{p\over q}$, and ${\bf O}'_{p\over q}$ respectively, where
$f_{A,b}({\bf p}_j)={\bf p}_{(j+1)_Q}$, $f_{A,b}({\bf p}'_j)={\bf
p}'_{(j+1)_Q}$, and ${P\over Q}=({p\over q})_1$.

  Assume that $b>1$. It follows from the
properties of ${\bf o}_{P\over Q}$ that
  the two critical points ${\bf
c}$ and ${\bf k}$ of $f_{A,b}$ are in an arc
  $\Gamma$ bounded by two
successive points ${\bf p}_j$ and ${\bf p}_k$ of ${\bf
  o}_{P\over
Q}$. When $Q=1$, ${\bf p}_j$ and ${\bf p}_k$ coincide. When
  ${\bf
o'}_{P\over Q}$ exists, they lie in an arc $\Gamma '$ bounded by two
successive
  points ${\bf p}'_j$ and ${\bf p}'_k$ of ${\bf o'}_{P\over
Q}$. Let $P_j$ be a lift
  of ${\bf p}_j$, $P_k$ be the lift of ${\bf
p}_k$ immediately to the right of $P_j$,
  and let $C$ and $K$ be the
lifts of ${\bf c}$ and ${\bf k}$ in $[P_j,\,P_k]$. Let
  $P'_j$ be the
lift of ${\bf p'}_j$ immediately to the left of $C$, and $P'_k$
be
  the lift of ${\bf p}'_k$ immediately to the right of $P'_j$ (and
of $K$). Let
  $F_{a,b}$ be the lift of $f_{A,b}$ such that $P_j$ and
$P'_j$ have rotation number
  ${p\over q}$ under $F_{a,b}$,
i.e.,
  $\underline{\rho}_F(P'_j)=\overline{\rho}_F(P'_j)={p\over q}$.
We then have the
  following result whose first part follows easily
from Theorem 3.2.1, and whose
  second part is a standard bifurcation
theory result.

  \medskip
\proclaim Lemma 3.6.3
{\rm (\cite{Bo,MaT})}. \\
(i)  With the above notation,
  $$(a,b)\in {\bf B}_{p\over
q}^l\setminus {\bf A}_{p\over
q}^l\Longleftrightarrow
  F_{a,b}(K)=F_{a,b}(P_j)\,,$$
  and
  $$(a,b)\in {\bf B}_{p\over q}^r\setminus {\bf A}_{p\over
q}^r\Longleftrightarrow
  F_{a,b}(C)=F_{a,b}(P_k)\,.$$
 (ii)  Furthermore,
for $b\not= 0$,
  $$(a,b)\in {\bf A}_{p\over q}^l\cap {\bf A}_{p\over
q}^r\Longleftrightarrow
   {\bf O}_{p\over
q}\,\mbox{is parabolic and }\,{\bf O}'_{p\over
  q}\,\mbox{ does not exist.}$$

  \medskip

  From the picture which emerges from the discussion so far, Theorem A
would follow from the uniqueness of the intersections described in
corollary 3.5.2.  The general case then follows by Theorem 3.3.1-3. We
shall prove this uniqueness property in \S~4 using the fact that the labels
of the curves that intersect determine the topological conjugacy classes of
the maps at the intersections.

\medskip
\proclaim  Lemma 3.6.4.
 At any crossing of two boundary curves ${\bf C}_{p\over
  q}^l$
and ${\bf D}_{p'\over q'}^r$, (where ${\bf C}$ and ${\bf D}$ stand for  either
${\bf
  A}$ or ${\bf B}$), the way the orbits ${\bf O}_{p\over q}$, ${\bf O}_{p'\over
  q'}$ and the critical points intertwine is
determined by the pair $({p\over
  q},{p'\over q'})$. Furthermore, the
itineraries of ${\bf O}_{p\over q}$ and ${\bf
  O}_{p'\over q'}$ are
determined by the pair $({p\over q},{p'\over q'})$.

  \medskip

\begin{figure}[ht]
	\centerline{\psfig{figure=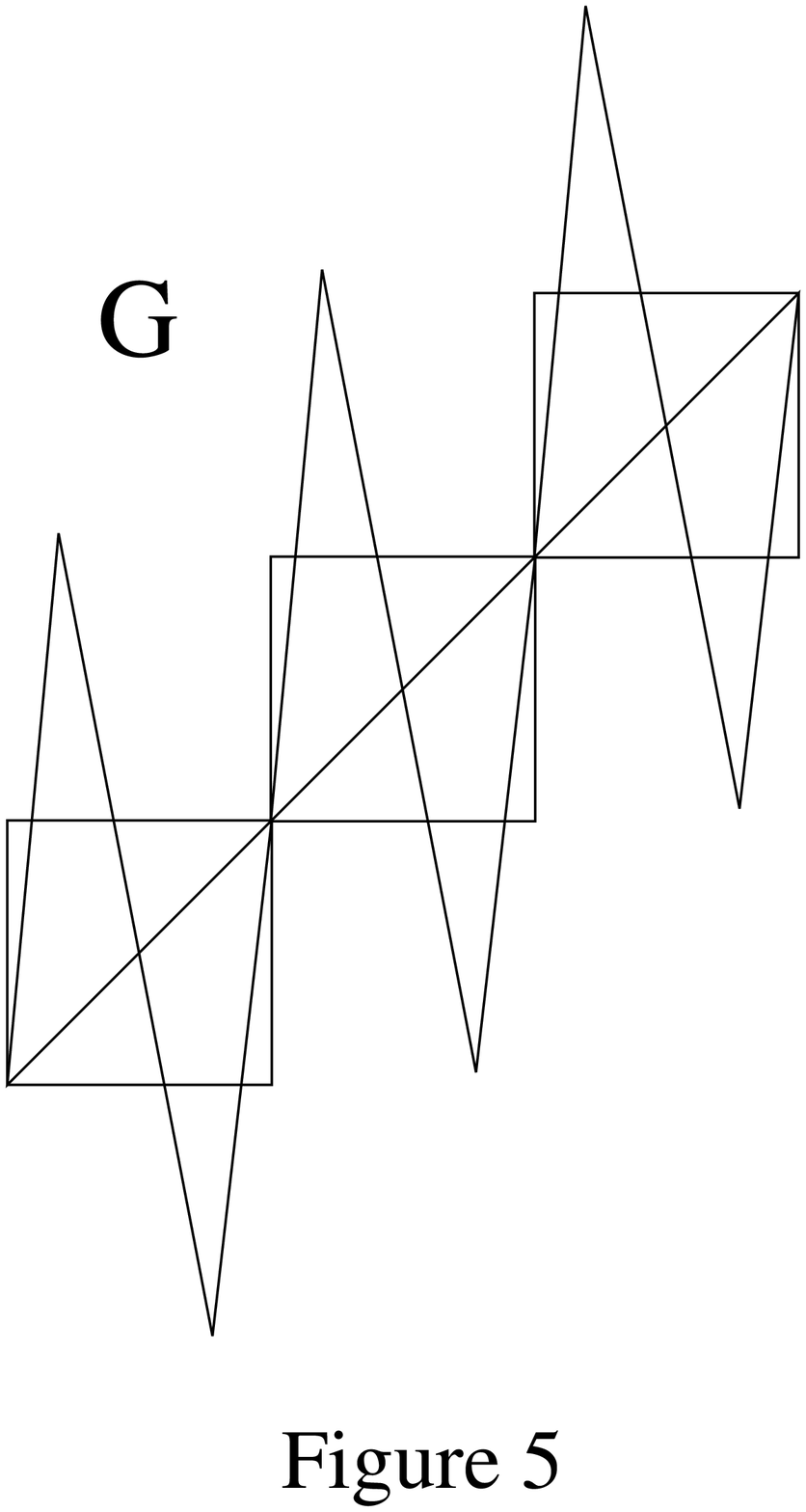,width=.3\hsize}}
\end{figure}

\proof  Take any two
standard lifts $F_{a_0,b_0}$ and $F_{a_1,b_1}$, which both possess the pair
of orbits $({\bf O}_{p\over q},{\bf O}_{p'\over q'})$. One can find a
piecewise linear lift $G$ as in Figure~5 that contains all periodic
itineraries of both $F_{a_0,b_0}$ and $F_{a_1,b_1}$.  Choosing $G$ to have
all its slopes greater than one in absolute value, it is easy to check that
for $\omega \in \lbrace {p\over{q}},{ p'\over{q'}}\rbrace$ it possesses only one
orbit which

 - Is invariant by a
non-decreasing lift with rotation number $\omega$,

- Has no point in the segments where $G$ is decreasing, and

- Is a lift of a periodic orbit of the circle map $g$.

 By standard kneading theory arguments we get(\cite{AM,MT}):

- The two orbits (for $\omega = {p\over q}$ and $\omega = {p'\over q'}$)
obtained this way and the turning points of $G$ are intertwined in the same
way as the corresponding orbits and critical points of $F_{a_0,b_0}$ and
$F_{a_1,b_1}$

- The kneading information about these orbits can be
read from $G$ as
  well.
\qed

  \medskip
\proclaim  Lemma 3.6.5. The maps corresponding to
 all intersections of the two boundary curves ${\bf C}_{p\over q}^l$ and
${\bf D}_{p'\over q'}^r$, (where again ${\bf C}$ and ${\bf D}$ stand for
either ${\bf A}$ or ${\bf B})$ are topologically conjugate.

  \medskip

\proof  This statement is a standard result of the topological
classification of maps with negative Schwarzian derivative,
  and
we refer to (\cite{MS} Chap. 2.3) for a more general discussion;
we give only a sketch of the arguments.

   Using Lemmas 3.6.3 and 3.6.4, we know that all such maps have the same
kneading data.  Because these are smooth maps with isolated critical
points, it follows that they have the same sets of itineraries. The fact
that maps with negative Schwarzian derivative off the critical set have no
homterval (intervals of positive length, not in the basin of a stable
periodic orbit, but where all iterates of the map are homeomorphisms
\cite{Ly,MMS}) yields
the conjugacy. Points with similar itineraries are paired by
  the
conjugacy, except for points belonging to the basins of stable or
semi-stable
  periodic orbits. For these points, the connected
components of the basins, on
  each side of the periodic points and
their preimages, can be paired in any way that
  respects the orbit
structure. Hence the conjugacy is not necessarily unique.
\qed

\section{Proof of Theorem A, Part II: Complex analytic part}

To complete the proof of Theorem A we must show that the boundary curves
described in \S~3.5 have unique intersection points; that is, that the
conjugacy classes in Lemma~3.6.5 correspond to a single map. This is the
content of Theorem D.  Before we can prove this lemma, however, we need to
introduce some techniques from complex analysis. References for the basic
theory of complex dynamics are
\cite{Kn7,Milnors-notes}. References for
Teichm\"uller theory are \cite{Gardiner,Lehto} and references for its
application to dynamical systems are
\cite{Epstein,Kn1,McM6,Sullivan}.

\subsection{Basic theory of complex dynamics}
We define a point to be {\it normal} for a family of holomorphic functions
if the functions in the family are locally uniformly bounded in a
neighborhood of the point. The set of normal points is open by definition.
We are interested in the normal sets of families generated by iterating a
single   holomorphic self-map of the punctured
plane $\Cstar$.

A {\em singular value} for a holomorphic map is either a {\em critical
value} (the image of a critical point) or an {\em asymptotic value} (a
limiting value of the image of a path tending to infinity).  A map with
only finitely many singular values is called a {\em finite type map}.  The
points $0$ and $\infty$ are asymptotic values for holomorphic self-maps of
$\Cstar$ but it will be more convenient to exclude them from the singular
value set.

The non-normal set divides the normal set into connected components. The
normal set is forward and backward invariant and its components are mapped
to one another. If a finite type map has a periodic cycle
$o=\{z_0,z_1=f(z_0), \ldots, z_{q-1}=f(z_{q-2})\}$ we define the multiplier
in the same way we do for real maps, that is, $m_o=f'(z_{q-2})\cdot
f'(z_{q-2}) \cdots f'(z_{0})$.

The following  theorem classifies the behavior of the components of the
normal set.

\proclaim Classification Theorem. Given a finite type
 holomorphic self map of $\Cstar$,
 the orbits of the components of the normal
set are characterized as follows:
\begin{itemize}
\item they fall onto
a  periodic cycle of components containg a periodic cycle with multiplier
$|m_o|< 1$ (attracting domain if $|m_o| \neq 0$
or super-attracting
domain if $|m_o| = 0$);
\item they fall onto  a periodic cycle
with
multiplier a root of unity (parabolic domain);
\item orbits
eventually
fall into a domain on which an iterate of the map is
holomorphically conjugate
to an irrational rotation (rotation domain).
 \end{itemize}

\Remark The classification of periodic normal behavior was
done by Fatou \cite{Fatou1,Fatou2}, Siegel \cite{Siegel} and Herman
\cite{Herman}.  The eventual periodicity of all normal components
(often called the Non-Wandering Theorem) was proved for rational maps by Sulliva
\cite{Sullivan}. For finite type holomorphic self-maps of
$\Cstar$ the Non-Wandering Theorem was proved in \cite{Kn1}.

Although arbitrary holomorphic self-maps of $\Cstar$ may have normal
components whose orbits fall onto a periodic cycle of domains in which
points are attracted to zero or infinity (essentially parabolic domains),
it was proved in
\cite{Kotus} that finite type holomorphic self-maps of $\Cstar$ have no
essentially parabolic domains.

\medskip
\Remark Each cycle  of periodic components  uses a
singular value in the following sense: cycles of super-attracting periodic
normal domains contain singular values by definition, cycles of attracting and
parabolic domains each contain the infinite forward orbit of a singular
value, and in fact one of the domains in the cycle contains the singular
point; finally, the boundary of any rotation domain is contained in the
closure of the forward orbit of some singular value. Proofs of these facts
go back to Fatou.  Among these facts is the statement in lemma 3.6.2-4(c).

\proclaim Definition.
The closure of the forward
orbits of the
singular values is called the {\em post-singular set} and is denoted by
$PS(f)$.

We shall be interested
in a special subclass of finite type maps.

\proclaim Definition.  A finite type
map is {\em geometrically finite} if
every infinite forward orbit of a singular value tends to a periodic cycle.

It may happen that no singular value has an infinite forward orbit; such
orbits are periodic or pre-periodic. These maps are trivially geometrically
finite.

Standard arguments (see e.g.
\cite{Epstein,Kn7,Milnors-notes}) show that for
geometrically finite maps

- There are no rotation domains, and

- Every infinite forward singular orbit
lies in the normal set and is attracted to a (necessarily attracting or
parabolic)  periodic cycle.

\subsection{Combinatorial Equivalence}

\proclaim Definition.
 A {\em combinatorial equivalence} of finite type maps is a pair of
homeomorphisms $(\phi,\psi)$ such that $$\phi \circ f_0 = f_1 \circ \psi$$
and such that $\phi$ and $\psi$ are isotopic rel $(PS(f_0))$.

\proclaim Definition. A homeomorphism $\phi:\Chat \rightarrow \Chat$ is called
{\em K-quasiconformal}, or {\em K-QC} for short, if there exists a $K \geq
1$ such that the field of infinitesimally small circles is mapped almost
everywhere onto a field of infinitesimally small ellipses of eccentricity
bounded by $k={K-1\over K+1}$. A map that is K-QC for some K is called {\em
quasiconformal}.

\proclaim Definition. A
combinatorial equivalence $(\phi,\psi)$ is {\em K-QC} (or just {\em QC} if
we don't care about the constant) if $\phi$ and $\psi$ are K-quasiconformal;
it is {\em strong} if $\phi$ and $\psi$ {\em agree} in a neighborhood of
each super-attracting, attracting and parabolic cycle (and hence define a
conjugacy in these neighborhoods).

Next we
prove,

\proclaim Lemma 4.2.1. A strong combinatorial equivalence  of
holomorphic geometrically finite maps can be isotoped (rel the post
singular set) through strong combinatorial equivalences to a strong QC
combinatorial equivalence.

\proof Let $(\phi,\psi)$ be the given
strong combinatorial equivalence. Let $N$ be the union of the neighborhoods
of the super-attracting, attracting and parabolic periodic cycles of $f_0$
on which $\phi$ and $\psi$ agree.

The first step is to isotop $\phi|_N=\psi|_N$ in $N$ rel $(PS(f_0) \cap N) \cup
\partial N$ to a quasiconformal homeomorphism that we again call $\phi$. To
do this, we use the
 canonical local picture  associated to  each
periodic cycle and determined by the multiplier of the cycle. For  a more
complete description of the local behavior  see e.g. \cite{Milnors-notes}.

\medskip

Case 1. Suppose first that $p$ is an attracting periodic point of $f_0$
with multiplier $\lambda$ and $N_p$ is the component of $N$ containing $p$.
Then there is an integer $k$ such that $f_0^k$ is the first return map for
$N_p$ and a conformal homeomorphism $h\colon N_p
\rightarrow \Delta$ where $\Delta$ is the unit disk, such that $h(0)=0$ and
$h \circ f_0^k(z) = \lambda h(z)$.

We can use the first return map $f_0^k$ to identify points in $N_p -
\lbrace p \rbrace$ and obtain a torus of modulus $\lambda$. The projection
$N_p - \{p\} \rightarrow N_p - \{p\}/f_0^k$ is a branched covering map of
the torus and the conjugation $\phi$ projects to this torus. Isotopy
classes rel the finitely many marked points for a torus are known to
contain K-QC maps for some $K > 1$, so the projection of $\phi$ may be
isotoped rel the branch points to a K-QC map. Since the homotopy lifting
property holds, (see e.g. \cite{GK1}), and the projection is holomorphic,
the K-QC map lifts to a K-QC map on $N_p - \{p\}$ and may be extended to
$N_p \cup \partial N_p$ so that the lift is isotopic to $\phi$ rel $(PS(f_0)
\cap N_p) \cup \partial N_p$.  (If, as in our
application to the standard family, there is a single branch point, and the
tori have the same modulus, the isotopy class contains a conformal map but
we do not use this fact.)

\medskip

Case 2. If $p$ is super-attracting the first return map is holomorphically
conjugate in $N_p$ to a map of the form $z \mapsto z^k$ on $\Delta$, for
some $k\geq 2$. If $p$ attracts no other singular points, we may push
$\phi$ to $\Delta$, isotop the map on $\Delta$ to a conformal map keeping
the boundary values fixed, and pull the isotoped map back to a conformal
map on $N_p$ rel $\partial N_p$. If $p$ does attract singular values the
argument has to be modified somewhat to take these orbits into account and
the isotopy will be only quasiconformal.  In our application we have two
singular values but we assume each is attracted to a {\em distinct}
periodic orbit. Hence we omit the details for the case where a
superattractive $p$ attracts a second singular value and refer the
interested reader to \cite{SulQCIII}.

\medskip

Case 3. It remains to describe the local behavior when $p$ belongs to a
parabolic cycle. The picture in this case is known as the Leau-Fatou
flower.  We make two simplifying assumptions: first that $p$ is a
non-degenerate parabolic fixed point, that is, $f_0'(p)=1$ and $f_0''(p)
\neq 0$, and second that $p$ attracts only one singular value. Full details
of the Leau-Fatou flower in the context of rational maps may be found in
\cite{Milnors-notes}, \S7. The details for geometrically finite maps may be
found in \cite{Epstein}.  A small neighborhood $N_p$ of $p$ is covered by a
pair of overlapping attracting and repelling {\em petals}, $U$ and
$U^{\prime}$, such that $f_0(U) \subset (U)$ and $f_0(U^{\prime}) \supset
U^{\prime}$.  If we conjugate $f_0$ by $w=-1/(z-p)$, $p$ is mapped to
infinity and the petals $U$ and $U^{\prime}$ are transformed into the two
overlapping regions $D_R$ and $D_L$ shown in figure 6. The conjugated map
in a neighborhood of infinity takes the form $$F
\colon w \mapsto w + 1 + o(1).$$
We see therefore that $D_L$ contains a left half plane $\{\Re w < - M \}$
for some $M>0$, in which $F$ is holomorphically conjugate to right
translation by $1$; similarly, $D_R$ contains a right half plane $\{\Re w >
M' \}$ for some $M'>0$ in which $F^{-1}$ is holomorphically conjugate to
left translation by $1$.

\begin{figure}[ht]
	\centerline{\psfig{figure=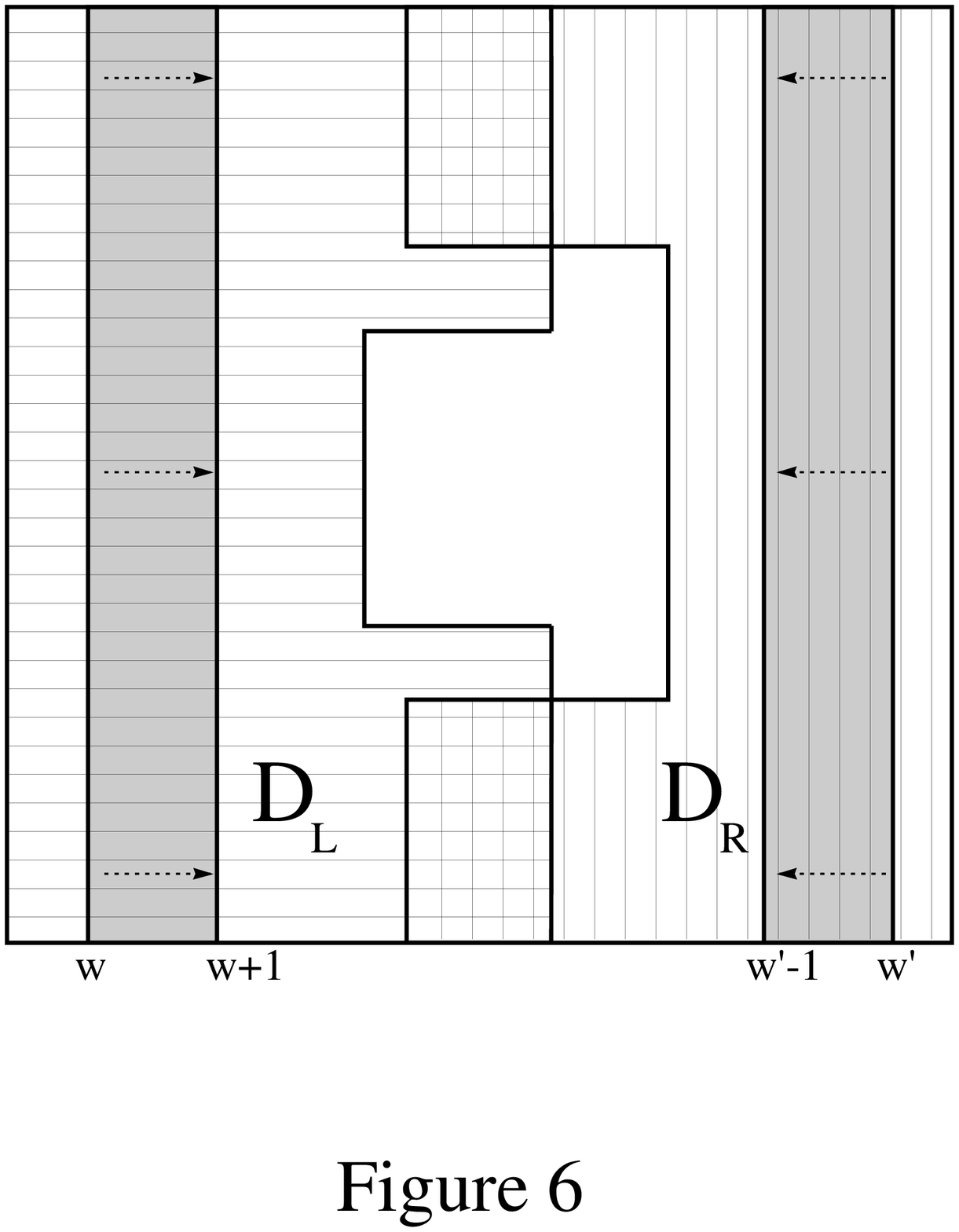,width=.4\hsize}}
\end{figure}

The orbit of the singular value in $U$ is transformed into the the
attracting region $D_R$ and since the map acts almost as translation the
imaginary parts of points in the orbit are bounded. The repelling petal
$U'$ is transformed into to the domain $D_L$. Hence $D_L$ contains a piece
of the non-normal set and so is not invariant under the conjugated map.

For each map $f_i, i=1,2$ we form the {\em Ecalle cylinder} $E_R$ by
identifying orbits in $D_R$ under the map $F$. The singular orbit projects
to a single marked point.  Similarly we form $E_L$ by identifying orbits in
$D_L$ under the map $F^{-1}$.  The conjugacy $\phi$ projects to the
cylinders and we can find, for some $K>1$, a K-QC map in the isotopy class
of this projected $\phi$.  Now we lift this isotoped $\phi$ to the exterior
of a large rectangle in the $w$-plane.

We thus have a quasiconformal conjugacy in a neighborhood of the parabolic
point (perhaps smaller than $N_p$). To obtain the strong K-QC equivalence
we must extend the conjugacy to a full neighborhood of the full
post-singular set.  Since $\phi$ is K-QC on $N$, where by definition $\phi
= f_1 \circ \phi \circ f_0^{-1}$, we may lift it as a K-QC map to
$f_0^{-1}(N)$. To extend this lift quasiconformally to the closures of
these neighborhoods we need to know that any
intersections of $\partial N$ and $\partial f_0^{-1}(N)$ are transverse. We
can assure this by modifying our original choice of $N$ if necessary, and
using the normal form for parabolic points again.
Since $f_0$ is geometrically finite, we may lift a finite number of times
to obtain a K-QC conjugacy on a neighborhood $N'$ of the full post-singular
set which is isotopic (rel $PS(f_0)$) to the original combinatorial
equivalence on $N'$ and agrees with it in the complement of $N'$.

\medskip

To complete the proof, we isotop $\phi$ in the complement of $N'$ to any
globally K-QC map and set $\psi = f_1^{-1} \circ \phi
\circ f_0$ where we choose the branch of the inverse to preserve the isotopy.
Note that these branches are well defined since there are no singular values in
this region.
\qed

\subsection{Application to the standard family}

  The circle maps $f_{A,b}$ have a natural extension to $\Cstar$. To see
this note that
the family of lifts $F_{a,b}$ extends to
 $\CC$, by the formula
 $$F_{a,b}:z \mapsto z + a + (b/2\pi)\sin(2\pi
z).$$   Using the projection of $\CC$ to $\Cstar$
given by the exponential map, we obtain
 holomorphic self-maps of $\Cstar$ that are holomorphic
extensions of the family $f_{A,b}$. For readability, we
keep the same notation. These maps
have exactly two
critical values and no asymptotic values so are of finite type.

\subsection{Extending real conjugacies}

The following lemma appears in various guises in
the literature. To prove a
version suited to our needs we require

\proclaim Definition. Let $I$ be an open interval in $\RR$ or $\TT$. A
homeomorphism $\phi:I \rightarrow I$ is called {\em K-quasisymmetric}, or
{\em K-QS}, if there exists a $K > 1$ such that for every triple $(a,b,c)$
of points in $I$, where $a < b < c$, $\phi$ satisfies $$ \frac{1}{K} <
\frac{\phi(c)-\phi(b)}{\phi(b) - \phi(a)} < K. $$ The homeomorphism $\psi:
\TT \rightarrow \TT$ is K-quasisymmetric if its restriction to every
subinterval is K-quasisymmetric.

\Remark The restriction of a K-QC homeomorphism of $\Cstar$ is
K-QS on $\TT$ and any K-QS homeomorphism of $\TT$
has a (not necessarily unique) K-QC extension to $\Cstar$ (see
\cite{Gardiner,Lehto}).

\proclaim Lemma 4.3.1.
Let $g_0,g_1$ be topologically conjugate maps of  $\TT$ in the family
$f_{A,b}$ whose extensions
$f_0,f_1$  to $\Cstar$ have the property that their post-singular
orbits remain in  $\TT$.  Then there is a strong K-QC combinatorial
equivalence $(\phi,\psi)$ for $f_0,f_1$.

\noindent\proof
We need only show  how to use the given real conjugacy $\Phi$ for
$g_0,g_1$ to obtain a strong combinatorial equivalence for $f_0, f_1$
because we may then apply Lemma 4.2.1 to complete the proof.

The first step is to replace $\Phi$ by a K-QS homeomorphism which
agrees with $\Phi$ on the closed post-singular set $PS(g_0)=PS(f_0)$.  We
can do this since $PS(g_0)$ consists of isolated points plus points
accumulating at attracting or parabolic cycles. The attracting and
parabolic cycles are distinguishable by their local topological behavior.
Near each cycle we use the local normal form to replace $\Phi$ by a
K-QS homeomorphism for some $K$; we then use the circle map $g_0$ to pull this
K-QS homeomorphism back to the closures of the basins of the
cycles in $\TT$; finally, we extend by continuity to $\TT$. Since $g_0$ is the
restriction of a holomorphic map, the new map, which we again call
$\Phi$, is K-quasisymmetric.

The second step is to extend the K-QS map $\Phi$ to a K-QC
 self-map $\phi$ of $\Cstar$. For each attracting or parabolic cycle
$P$, let $N_P$ be a neighborhood of $P$ in $\Cstar$ with smooth boundary.
Using the local normal form again, we extend the K-QS map $\Phi$ to
$N_P$ so that it is K-QC. Extending this way for all the cycles
defines a germ $\phi$ for a K-QC conjugacy between $f_0$ and $f_1$. Now we
extend $\phi$ arbitrarily as a K-QC homeomorphism of $\Cstar$.

The final step is to define a lift $\psi$ of $\phi$ so that the pair
$(\phi,\psi)$ are isotopic rel the post-singular sets and are the desired
strong combinatorial equivalence. Denote the critical value set of the map
$f_i$ by $S_i$, $i=0,1$; each set consists of two points.  The maps $f_i$
are covering maps of $\Cstar - f_i^{-1}(S_i)$ onto $\Cstar - S_i$. Extend
these covering maps to fix the ``ends'' zero and infinity of $\Cstar$.
Because the maps $f_0,f_1$ are in the same family, that is, given by a
formula of the form $\alpha
\zeta \exp \beta(\zeta - 1/ \zeta)$ for constants $\alpha$ and $\beta$,
and variable $\zeta \in \Cstar$, they are built up from a sequence of
elementary maps whose lifting properties are known.  The lift $\psi =
f_1^{-1} \circ
\phi \circ  f_0$ may therefore be defined uniquely so that it
 agrees with $\phi$ on any and hence all the points of $PS(f_0)$. \qed

\medskip

\Remark If the post-singular set is actually finite, the situation
is much simpler. Every singular point is superattracting or else its orbit
eventually lands on a repelling periodic cycle. We can choose an arbitrary
topological extension to $\Cstar$ as the homeomorphism $\phi$ and define
$\psi$ by the formula $\psi = f_1^{-1} \circ \phi \circ f_0$
where again the branch of the inverse is chosen so that $(\phi,\psi)$ are
isotopic rel the post-singular sets.  By the easy parts of Lemma 4.2.1
there are automatically quasiconformal homeomorphisms in this isotopy
class.

\subsection{Statement of the Rigidity Theorem}
\label{sec:statement}

\proclaim Theorem D (Rigidity). Suppose that the functions $f_0$
and $f_1$ are both intersections of boundary curves  $C^l_{p\over q}$ and
$D^r_{p' \over q'}$ where $C$ and $D$ stand for either $A$ or $B$ as in Lemma
3.6.5. Then $f_0=f_1$.

>From Lemma 3.6.3 we see that at the intersections of the boundary curves
one of the following holds:

 $\alpha$ - Both singular orbits are attracted by distinct parabolic cycles,

 $\beta$ - Both singular orbits are preperiodic, or

 $\gamma$ - One singular orbit is preperiodic and the other is attracted by
a parabolic cycle.

It follows that the extensions of standard maps corresponding to these
intersection points are geometrically finite.

\subsection{Basic Teichm\"uller theory}

To prove the Rigidity Theorem D we
follow the version of the proof of a rigidity result for rational maps carried
out by McMullen in \cite{McM6}. In particular, we shall use some standard
Teichm\"uller theory.  Below we state those facts we require  in a
form suited to our needs. A good basic reference for this material is
\cite{Gardiner}, Chap. 6.  Thurston and
Sullivan were the first to apply these techniques in the context of
rigidity in holomorphic dynamics.

\medskip

Let $X$ be a compact Riemann surface and let $C$ be a closed subset of $X$
containing at most countably many points.
  Then
  the {\em Teichm\"uller space} of  $X$ with boundary $C$
is the set of isotopy classes of quasiconformal homeomorphisms of $X$ rel
$C$. We denote it by $\TTT(X,C)$.  We shall be interested in $\TTT(X,C)$
where $X =\Cstar$ and $C= PS(f)$ for $f$ in the
standard family.

The Teichm\"uller space is finite dimensional if $X$ has finite genus and
$C$ is a finite point set:
 in our case, if $PS(f)$ is finite.

\medskip

If $\phi$ is a quasiconformal homeomorphism of $X$, its {\em Beltrami
differential}
is $\mu(z) =
\phi_{\bar{z}}/\phi_{z}$, where the derivatives are taken in the
generalized sense.  The infinitesimal ellipse field is determined by
$\mu(z)$: the eccentricity of the ellipse at the point $z$ is $|\mu(z)|$
and the major axis  has argument $\arg \mu(z)$.

The {\em maximal dilatation} of $\phi$ is
$$K_{\phi}(X)= \max_z (1 + ||\mu
||_{\infty})/ (1 - ||\mu ||_{\infty}) < \infty. $$

Given an isotopy class  of quasiconformal homeomorphisms $X$ rel $C$
 one can
ask  if there is a map that is {\em extremal}; that is, its maximal
dilatation is minimal over all maps in its class.

\medskip

 A quasiconformal map is called a {\em Teichm\"uller map} if it is locally an
affine stretch: that is, its Beltrami differential  has the form $\mu = t
{\bar q}/|q|$ where $q$ is a holomorphic quadratic differential such that
$||q||=\int_X |q| <
\infty$ and $|t| < 1$.

\proclaim Teichm\"uller's Theorem. Let $\TTT(X,C)$ be a finite dimensional
Teichm\"uller space. Then every isotopy class contains an extremal map.
Moreover, this extremal is unique and is a Teichm\"uller map. If $\TTT(X,C)$
is not finite dimensional, the extremal map exists but it is not
necessarily unique nor is any such extremal a Teichm\"uller map.

Since the post-singular set is not always finite we need to consider
infinite dimensional Teichm\"uller spaces.  To this end, we introduce the
concept of boundary dilatation. Let $S = X-C$ and let $R$ be any compact
subset of $S$. Set $K_{\phi}^0(S-R) = \inf_{\psi \sim
\phi}(K_{\psi}(S-R))$.  Define the {\it boundary dilatation} $H(\phi)$ as
the direct limit of the numbers $K_{\phi}^0(S-R)$ as $R$ increases to $S$.

\proclaim Strebel's Frame Mapping Condition. Let $\phi$ be a quasiconformal
homeomorphism of $S$ to another surface and suppose $H(\phi) <
K^0_{\phi}(S)$. Then the isotopy class of $\phi$ (rel $C$) contains a
unique extremal map and this map is a Teichm\"uller map.

\subsection{Proof of Theorem D}

 It suffices to prove that if $f_0$ and $f_1$ are topologically conjugate
maps in the standard family whose singular orbits satisfy one of the
conditions $\alpha - \gamma$ of section
\ref{sec:statement} then they are equal.

 Since their extensions to $\Cstar$ are geometrically finite, by Lemma
4.3.1 there is a strong K-QC combinatorial equivalence $(\phi,\psi)$
between them.

 Suppose first that both singular orbits are preperiodic.  Then the
post-singular set is finite and any K-QC combinatorial equivalence is
trivially strong.  Moreover, $\TTT(\Cstar,PS(f_0))$ is finite dimensional
and by Teichm\"uller's theorem, there is a unique extremal map in every
isotopy class; denote the extremal map in the isotopy class of $\phi$ and
$\psi$ by $\hat\phi$.  Now we replace $\phi$ by $\hat\phi$ as we did in the
last step of the proof of Lemma 4.2.1 and set $\hat\phi = f_1^{-1} \circ
\hat\phi
\circ f_0$,  choosing the branch that preserves the isotopy. Since $f_0$ and
$f_1$ are holomorphic the infinitesimal ellipse fields determined by
$\hat\phi$ and $\hat\psi$ are the same and $\hat\psi$ is extremal. By
uniqueness $\hat\phi=\hat\psi$;  denote the extremal quasiconformal
conjugacy $\hat\phi$ by $\phi$ again.

\medskip

In the other two cases, there is at least one singular orbit attracted by a
parabolic cycle. It is important to note that no parabolic cycle attracts
more that one singular orbit.  The quasiconformal homeomorphisms $(\phi,
\psi)$ we obtained in the proof of Lemma 4.2.1 agree in a neighborhood $N$ of
the post-singular set.  We need to modify this $\phi$ in a neighborhood of
a parabolic point $p$ containing the forward orbit of one singular value so
that it satisfies the Frame Mapping Condition.

 As above we conjugate $f_0$ to $F(w)=w + 1 + o(1)$ by sending $p$ to
infinity.  We follow the argument in \cite{Epstein}, \S 4.2, lemma 78.  An
application of the Schwarz lemma shows that $|F'(w)|$ is uniformly close to
$1$ in a neighborhood of infinity.  This means that for $\eta$ large, the
image of $F(t \pm i \eta)$, $ t \in \RR$ is a curve that stays very close
to horizontal. Hence, given any $\epsilon > 0$, we can find $M$ such that
for $|\eta| > M$, $\phi$ is isotopic to a map (again called $\phi$) with
dilatation less than $1 +
\epsilon$. Next, using the images of the endpoints of vertical lines inside
the closed large rectangle to control the images of these lines, and noting
that we have arranged it so that there are no  points of $PS(f)$ inside the
large rectangle, we can isotop $\phi$ in the part of $D_L
\cup D_R$  inside the rectangle so that it is
quasiconformal.

This new map together with its lift in the same isotopy class gives us a
combinatorial equivalence $(\phi,\psi)$ which is no longer strong but is
still K-QC.  This new $\phi$ satisfies Strebel's Frame Mapping Condition
for $\epsilon$ small enough. Therefore, just as in the preperiodic case, we
may replace both maps in the equivalence with the unique extremal
Teichm\"uller map in their isotopy class and obtain a quasiconformal
conjugation, denoted again by $\phi$.

\medskip

Finally, we complete the proof of the lemma by showing that $\phi$ is
conformal and hence a homothety.

If $\phi$ is not conformal, its Beltrami differential determines a
quadratic differential $q$ on $S $.  Since $\phi$ is a conjugacy, and the
maps $f_0$ and $f_1$ are holomorphic, the infinitesimal ellipse fields
determined by the Beltrami differential $\mu$ and the
Beltrami differential $f_0^* \mu$ of $f_1^{-1} \circ\, \phi\, \circ f_0 = \phi$
are the same; that is,
 $f_0^* \mu = \mu$.
Since $f^*_0 \mu$ is again the Beltrami differential of a Teichm\"uller
map, it has the form $ f_0^*\mu
= t \overline{f^*_0 q}/|f_0^* q|$  where  $f_0^*q$ is the  pull-back quadratic
differential.  Now on the one hand,  the norm of the pullback differential
$||f_0^* q||$ is given by $||q||$ times the degree of $f_0$ so since $f_0$
has infinite degree, $||f_0^* q||$
is unbounded. On the other hand however, $f^*_0
\mu = t \overline{f^*_0 q}/|f^*_0 q| $, so that $$\bar{f^*_0 q}/|f^*_0 q|=
\bar q/|q|.$$ If $h=f^*_0 q /q$, then $\bar h = |h|$ and $h$ is real valued.
But $h$ is meromorphic on $S$ and any meromorphic function taking only real
values  must be constant. Thus $f^*_0 q = c q$
for some $c>0$. Since $||q||$ is bounded, we have a contradiction and
$\phi$ is conformal.

\medskip

If we conjugate $f_{A,b}$ by a homothety, we obtain an equivalent dynamical
system.  Since the homothety preserves the unit circle, the factor must
have modulus $1$; its argument only appears in the sine term and doesn't
change any of the dynamical properties. \qed

\section{Concluding remarks}
In our study of the standard family we used
real analytic techniques  to get good control of the
boundary curves in the parameter plane of regions with a given lower or
upper bound on the rotation number. In order to control the intersections
of these curves  we needed to apply  rigidity properties found in
families of complex analytic maps.
 Previously,  complex analytic techniques were used to obtain
rigidity  in  one parameter
 families of maps with a single critical point.
 Our description of the
parameter space of the standard family is still incomplete and we pose some
open problems here. They do not  seem amenable to the methods
used so far and new ideas are needed.
\newpage

\noindent
{\bf Conjectures:}
\begin{itemize}
\item
$R_{p_0\over q_0}$ is homeomorphic to $R_{p_1\over q_1}$ by a homeomorphism
 $H_{{p_0\over q_0},{p_1\over q_1}}$ having the following property:\\
 If $H_{{p_0\over q_0},{p_1\over q_1}}(f_0)=f_1$, then $f_0$  partitions
the circle into $q_0$ intervals, $I_1,\ldots I_{_0}$, and $f_1$ partitions
the circle into $q_1$ intervals,
$J_1,\ldots J_{q_1}$ so that $f_0^{q_0}|_{I_j}$ is topologically
conjugate to $f_1^{q_1}|_{I_k}$, $j \in \{1,\ldots,q_0\}, k\in
\{1,\ldots,q_1\}$.

\item
The set of maps with a given topological entropy is connected.
\end{itemize}

Both conjectures are proved  in \cite{UTGC} for some two parameter families
of piecewise affine maps.  A similar entropy conjecture for cubic maps is
discussed in \cite{DGMT}.

\bibliography{refer}
 \end{document}